\documentclass[12pt]{article}
\usepackage{amssymb,latexsym}
\usepackage{fullpage}
\usepackage{array}
\renewcommand{\thefootnote}{}

\begin{document}
\title{K-theory of quasi-toric manifolds}
\author{P.Sankaran and V.Uma}
\date{}
\maketitle
\thispagestyle{empty}
\thefootnote{{\bf AMS Subject Classification:} 
Primary: 55N15, Secondary: 14M25\\
{\bf keywords:}  Quasi-toric manifolds, K-theory, 
Bott towers, Bott-Samelson varieties}

\def\theequation {\arabic{section}.\arabic{equation}}
\renewcommand{\thefootnote}{}

\newcommand{\codim}{\mbox{{\rm codim}$\,$}}
\newcommand{\stab}{\mbox{{\rm stab}$\,$}}
\newcommand{\lr}{\mbox{$\longrightarrow$}}

\newcommand{\ch}{{\cal H}}
\newcommand{\cf}{{\cal F}}
\newcommand{\cd}{{\cal D}}

\newcommand{\blr}{\Big \longrightarrow}
\newcommand{\da}{\Big \downarrow}
\newcommand{\ua}{\Big \uparrow}
\newcommand{\hra}{\mbox{\LARGE{$\hookrightarrow$}}}
\newcommand{\rt}{\mbox{\Large{$\rightarrowtail$}}}
\newcommand{\dua}{\begin{array}[t]{c}
\Big\uparrow \\ [-4mm]
\scriptscriptstyle \wedge \end{array}}

\newcommand{\be}{\begin{equation}}
\newcommand{\ee}{\end{equation}}

\newtheorem{guess}{Theorem}[section]
\newcommand{\bth}{\begin{guess}$\!\!\!${\bf .}~}
\newcommand{\eeth}{\end{guess}}
\renewcommand{\bar}{\overline}
\newtheorem{propo}[guess]{Proposition}
\newcommand{\bpropo}{\begin{propo}$\!\!\!${\bf .}~}
\newcommand{\epropo}{\end{propo}}

\newtheorem{lema}[guess]{Lemma}
\newcommand{\blem}{\begin{lema}$\!\!\!${\bf .}~}
\newcommand{\elem}{\end{lema}}

\newtheorem{defe}[guess]{Definition}
\newcommand{\bdefe}{\begin{defe}$\!\!\!${\bf .}~}
\newcommand{\edefe}{\end{defe}}

\newtheorem{coro}[guess]{Corollary}
\newcommand{\bcor}{\begin{coro}$\!\!\!${\bf .}~}
\newcommand{\ecor}{\end{coro}}

\newtheorem{rema}[guess]{Remark}
\newcommand{\brem}{\begin{rema}$\!\!\!${\bf .}~\rm}
\newcommand{\erem}{\end{rema}}

\newtheorem{exam}[guess]{Example}
\newcommand{\beg}{\begin{exam}$\!\!\!${\bf .}~\rm}
\newcommand{\eeg}{\end{exam}}

\newcommand{\ctext}[1]{\makebox(0,0){#1}}
\setlength{\unitlength}{0.1mm}
\newcommand{\cl}{{\cal L}}
\newcommand{\cp}{{\cal P}}
\newcommand{\ci}{{\cal I}}
\newcommand{\bz}{\mathbb{Z}}
\newcommand{\cs}{{\cal s}}
\newcommand{\cv}{{\cal V}}
\newcommand{\ce}{{\cal E}}
\newcommand{\ck}{{\cal K}}
\newcommand{\cR}{{\cal R}}
\newcommand{\cz}{{\cal Z}}
\newcommand{\cg}{{\cal G}}
\newcommand{\bq}{\mathbb{Q}}
\newcommand{\bt}{\mathbb{T}}
\newcommand{\bh}{\mathbb{H}}
\newcommand{\br}{\mathbb{R}}
\newcommand{\wt}{\widetilde}
\newcommand{\im}{{\rm Im}\,}
\newcommand{\bc}{\mathbb{C}}
\newcommand{\bp}{\mathbb{P}}
\newcommand{\spin}{{\rm Spin}\,}
\newcommand{\ds}{\displaystyle}
\newcommand{\tor}{{\rm Tor}\,}
\newcommand{\bff}{{\bf F}}
\newcommand{\bs}{\mathbb{S}}
\def\ns{\mathop{\lr}}
\def\nssup{\mathop{\lr\,sup}}
\def\nsinf{\mathop{\lr\,inf}}
\renewcommand{\phi}{\varphi}
\newcommand{\co}{{\cal O}}
\noindent
\begin{abstract} In this note we shall give 
a description of the $K$-ring of a quasi-toric manifold 
in terms of generators and relations. We apply 
our results to describe the $K$-ring of Bott-Samelson varieties. 
\end{abstract}

\section{Introduction}
The notion of quasi-toric manifolds  
is due to M.Davis and T.Januszkiewicz \cite{dj} who 
called them `toric manifolds'. The quasi-toric manifolds 
are a natural topological generalization  of the algebraic geometric 
notion of non-singular projective toric varieties.    
However 
there are compact complex non-projective non-singular toric varieties which 
are quasi-toric manifolds. See \cite{bp}.  Recently 
Civan \cite{c} has constructed 
an example of a compact complex non-singular toric variety 
which is {\it not} a quasi-toric manifold. 
 
In \cite{su}, we obtained, among other things,  
a description of the $K$-ring of projective non-singular 
toric varieties in terms of generators and relations. 
(In fact our result is applicable to slightly more general 
class of varieties.) The purpose of this note is to 
extend the $K$-theoretic results of \cite{su} to the context of 
quasi-toric manifolds.  As an application we obtain 
a description of the K-ring of Bott-Samelson varieties. 
The more difficult problem of computing  
the $KO$-theory has been solved by A.Bahri and M.Bendersky. 

Let $G=(\bs^1)^n$ be an $n$-dimensional 
{\it compact} torus and let $P\subset \br^n$ be a simple convex 
polytope of dimension $n$. That is, $P$ is a convex 
polytope in which exactly $n$ 
facets --- codimension $1$ faces of $P$ --- meet at each vertex of 
$P$. A $G$-quasi-toric manifold over $P$ 
is a (smooth) $G$-manifold $M$ where the $G$-action  
is locally standard with projection $\pi:M\lr M/G\cong P$.   
Here `local standardness' means that every point of $M$ has 
an equivariant neighbourhood $U$ such that there exists an 
automorphism $\theta:G\lr G$, an equivariant open subset $U'\subset \bc^n$ 
where $G$ action on $\bc^n$ is given by the standard inclusion 
of $G\subset U(n)$, and a diffeomorphism $f:U\lr U'$ 
where $f(t\cdot x)=\theta(t)f(x)$ for all $x\in U, t\in G$. 
Any two points of $\pi^{-1}(p)$ have the 
same isotropy group its dimension being  
codimension of the face of $P$ which contains $p$ in its relative 
interior. It is known that $M$ admits a CW-structure  
with only even dimensional cells. In particular 
$M$ is simply connected and hence orientable.  

Let $\cf_P$  (or simply $\cf$) denote 
the set of facets of $P$ and let $|\cf|=d$. For each 
$F_j\in \cf$, let $M_j=\pi^{-1}(F_j)$ and $G_j$ be the 
($1$-dimensional) isotropy subgroup at any `generic' point of $M_j$. 
Then $M_j$ is orientable for each $j$. 
The subgroup $G_j$ determines a primitive vector 
$v_j$ in $\bz^n=Hom(\bs^1,G)$ which 
is unique upto sign.  The sign is determined by 
choosing an omni-orientation on $M$, i.e. orientations on $M$ as well as 
one on each $M_j$, $1\leq j\leq d$. 
Choosing such a $v_j$ for $1\leq j\leq d$ 
defines the `characteristic map' $\lambda:\cf\lr \bz^n\cong Hom(\bs^1, G)$ 
where $F_j\mapsto v_j$.  Suppose that $F_1,\cdots, F_d$ are 
the facets of $P$, then writing $v_j=\lambda(F_j)$, the 
primitive vectors $v_1,\cdots, v_d$ are such that 
$$ \textrm{if }\bigcap_{1\leq r\leq k}F_{j_r}\subset P \quad 
\textrm{is of codimension $k$ then }$$
$$v_{j_1},\cdots, v_{j_k} 
\textrm{ extends to a $\bz$-basis 
$v_{j_1},\cdots, v_{j_k},w_1,\cdots, w_{n-k}$ of }\bz^n. 
\eqno(1.1)$$ 
Fix an orientation for $M$. The omni-orientation on $M$ determined by 
$\lambda$ is obtained by orienting $M_j$ so that the oriented normal 
bundle corresponds to the $1$-parameter subgroup given by $v_j.$

We shall call any map $\lambda:\cf_P\lr \bz^n$ that satisfies 
$(1.1)$ a characteristic map.   

Conversely, starting with a pair $(P, \lambda)$ 
where $P$ is any simple convex polytope and 
a characteristic map $\lambda:\cf\lr \bz^n$  
there exists a quasi-toric 
manifold $M$ over $P$ whose characteristic 
map is $\lambda$.  The data $(P,\lambda)$ determines 
the $G$-manifold $M$ and an omni-orientation on it.   
We refer the reader to \cite{dj} and \cite{bp} for basic  
facts  concerning quasi-toric manifolds. 

Suppose that $P$ is a simple convex polytope of dimension 
$n$ and that $F_j\mapsto v_j, 1\leq j\leq d$ is a 
characteristic map $\lambda:\cf\lr \bz^n$. Assume 
that $F_1\cap\cdots \cap F_n$ is a vertex of $P$ 
so that $v_1,\cdots,v_n$ is a $\bz$-basis of $\bz^n$. 
Let $S$ be any commutative ring with identity and let 
$r_1,\cdots, r_n$ be invertible in $S$. 

\bdefe \label{ring}
Consider 
the ideal $\ci$ of the polynomial algebra 
$S[x_1,\cdots,x_d]$ generated by the following 
two types of elements:
$$ x_{j_1}\cdots x_{j_k}\eqno{(1.2)}$$
whenever $F_{j_1}\cap \cdots \cap F_{j_k}=\emptyset, $  
and the elements 
$$z_u:=\prod_{j, u(v_j)>0} (1-x_j)^{u(v_j)}-r_u \prod_{j, u(v_j)<0}
(1-x_j)^{-u(v_j)}\eqno{(1.3)}$$ 
where $u\in Hom(\bz^n,\bz)\cong Hom(G,\bs^1)$ and 
$r_u:=\prod_{1\leq i\leq n}r_i^{u(v_i)}.$  We denote the 
quotient $S[x_1,\cdots,x_d]/\ci$ by $\cR(S;\lambda)$ 
or simply by $\cR$. 
\edefe

Let $E\lr B$ be a principal $G$-bundle with base space 
$B$ a compact Hausdorff space. Denote by $E(M)$ the 
associated $M$-bundle with projection map $p:E\times_GM\lr B$.  
The choice of the basis $v_1\cdots, v_n$ for $\bz^n=Hom(\bs^1,G)$ 
yields a product decomposition $G=\prod_{1\leq i\leq n}G_i\cong(\bs^1)^n$.  
Also one obtains principal $\bs^1$-bundles 
$\xi_i,~1\leq i\leq n,$  over $B$ associated 
to the $i$-th projection $G=(\bs^1)^n\lr \bs^1$. The projection    
$E\lr B$ is then the projection of the bundle 
$\xi_1\times \cdots \times\xi_n$ over $B$. For any $G$-equivariant 
vector bundle $V$ over $M$ denote by $\cv$ the bundle over $E(M)$ 
with projection $E(V)\lr E(M)$.  We shall often denote the complex line 
bundle associated to a principal $\bs^1$ 
bundle $\xi$ also by the same symbol $\xi$. 

Suppose that $V$ is the product 
bundle $M\times \bc_\chi$, where $\bc_\chi$ is the $1$-dimensional 
$G$-representation given by the character $\chi:G\lr \bs^1$. Then 
$\cv$ is isomorphic to the pull-back of the bundle $p^*(E_\chi)$, where 
$E_\chi$ is obtained from $E\lr B$ by 
`extending' the structure group to $\bs^1$ via the character $\chi$.  
Writing $\chi=\sum_{1\leq i\leq n} a_i\rho_i$, where $\rho_i:G\lr \bs^1$ 
is the $i$-th projection, one has $\cv\cong 
p^*(\xi_1^{a_1}\cdots \xi_n^{a_n})$, where $\xi^a=(\xi^*)^{-a}$ when $a<0$. 
 
\bth \label{main}
Let $M$ be a quasi-toric manifold 
over a simple convex polytope $P\subset \br^n$ and characteristic 
map $\lambda:\cf\lr \bz^n$.  Let $E\lr B$ be a 
principal $G=(\bs^1)^n$ bundle over a compact Hausdorff space.  
With the above notations,  
there exist equivariant line bundles $L_j$ over $M$ such 
that, setting $r_i=[\xi_i]\in K(B),~1\leq i\leq n$, 
one has an isomorphism 
of $K(B)$-algebras $\phi:\cR(K(B);\lambda)\lr K(E(M))$ defined by 
$x_j\mapsto (1-[\cl_j])$. 
\eeth 

The proof is given in \S 3. A technical lemma needed in 
the proof is established in \S2. In \S4, we apply 
our result to obtain the $K$-ring of   
Bott-Samelson varieties.

%%%%%%%%%%%%%%%%%%%%%%%%%%%%%%%%%%%%%%%%%%%%%%%%%%%%%%%%%%%
%%%%%%%%%%%%%%%%%%%%%%%%%%%%%%%%%%%%%%%%%%%%%%%%%%%%%%%%%%%
%%%%%%%%%%%%%%%%%%%%%%%%%%%%%%%%%%%%%%%%%%%%%%%%%%%%%%%%%%%

\section{Generators of $\cR$}

We keep the notations of the previous section. In this 
section we give a convenient generating set for $S$-module $\cR(S;\lambda)$ 
where $S$ is any commutative ring with identity and
$\lambda:\cf_P\lr \bz^n$ a characteristic map, $P\subset \br^n$ being a 
simple convex $n$ dimensional polytope. 

Let $h:\br^n\lr \br$ be a linear map which is injective 
when restricted to the set $P_0$ of vertices of $P$.
Then $h$ is a generic ``height function'' with 
respect to the polytope $P$. That is, $h$ is injective when 
restricted to any facet of $P$.  The map $h$ induces 
an ordering on the set $P_0$ of vertices of $P$, where $w<w'$ in $P_0$ if 
$h(w)<h(w')$.  The ordering on $P_0$ induces an orientation on the 
edges $P_1$ of $P$ in the obvious fashion.  Since $P$ is simple, 
there are exactly $n$ edges which meet at each 
vertex of $P$.  Given any $w\in P_0$, denote by $T_w$ the face 
of $P$ spanned by all those edges incident at $w$ which point {\it away} 
from $w$.   Then the following property holds: 
$$\textrm{if}\quad w'\in P_0 ~~\textrm{belongs to} \quad T_w, ~\textrm{then} 
\quad w\leq w'.  
\eqno(2.1)  $$ 
This is a consequence of the assumption that $P$ is simple and can be 
proved easily using for example Lemma 3.6 of \cite{z}. Property 
(2.1) is `dual' to  property $(*)$ of \S5.1, \cite{f}.  

It is shown in \cite{dj} that 
$M$ has a perfect cell decomposition with respect to 
which the submanifolds 
$M_w:=\pi^{-1}(T_w)$ are closed cells. (Cf. \cite{b}, and \cite{kh}.)

If $Q$ is a proper face of $P$, we denote by $x(Q)$ the product 
$x_{j_1}\cdots x_{j_r}$ where $F_{j_1}, \cdots, F_{j_r}$ are the 
distinct facets of $P$ which contain $Q$. 

\blem \label{gen}  
With the above notation, the elements $x(T_w),~w\in P_0,$ generate 
$\cR$ as an \newline $S$-module.  \hfill $\Box$   
\elem 
\noindent 
The proof is identical to that of lemma 2.2(iv), \cite{su} 
using (2.1) in the place of property $(*)$ of \cite{su}.  
Proof of lemma 2.2(i), \cite{su}, which was omitted, had an error.
(See Errata \cite{su}.)  However it is redundant 
here since $z_u=0$ in $\cR$ in view of equation (1.3) 
in Definition \ref{ring} above. 

%%%%%%%%%%%%%%%%%%%%%%%%%%%%%%%%%%%%%%%%%%%%%%%%%
%%%%%%%%%%%%%%%%%%%%%%%%%%%%%%%%%%%%%%%%%%%%%%%%%

\section{Proof of Theorem \ref{main}}
We keep the notations of the previous section.  Let $M$ be 
the quasi-toric manifold over a simple convex polytope 
$P\subset \br^n$ with characteristic map $\lambda:\cf\lr \bz^n$. 
Let $|\cf|=d$.  As in the previous section we 
shall assume that $\bigcap_{1\leq i\leq n}F_i$ is a vertex 
of $P$ and $G=\prod_{1\leq i\leq n}G_i\cong(\bs^1)^n$ the corresponding 
product decomposition.  

Set $\wt{G}=(\bs^1)^d$ and let $\theta:\cf\lr \bz^d$ be defined 
by  $F_j\mapsto e_j,$ the standard basis vector, for each $j\leq d$. 
For any face $F=F_{j_1}\cap \cdots \cap F_{j_k}$, set 
$\wt{G}_F$ denote the subgroup $\{(t_1,\cdots, t_d)\in \wt{G}\mid 
t_i=1,\quad i\neq j_1,\cdots, j_k\}$. 
One has a $\wt{G}$-manifold 
$\cz_P:=\wt{G}\times P/\sim$ where $(g,p)\sim (g',p')$ 
if and only if $p=p'$ and $g^{-1}g'\in \wt{G}_F$ where 
$p$ is in the relative interior of the face $F\subset P.$
The action of $\wt{G}$ on $\cz_P$ is given by $g.[g',p]
=[gg',p]$ for $g\in \wt{G}$ and $[g',p]\in \cz_P$. 
One has the projection map $p:\cz_P\lr P$.  However $\cz_P$
is not a quasi-toric manifold over $P$ since $d>n$.  When $P$ is 
clear from the context we shall denote $\cz_P$ simply by $\cz$. 

Let $\wt{\lambda}:\bz^d\lr \bz^n$ be defined 
as $e_j\mapsto v_j:=\lambda(F_j), 
\quad 1\leq j\leq d$.
This corresponds to a surjective 
homomorphism of groups $\Lambda:\wt{G} \lr G$ with kernel 
$H\subset \wt{G}$ for the subgroup  
corresponding to $\ker(\wt{\lambda})$.   
One has a splitting $\bz^d=\ker(\wt{\lambda})\oplus 
\bz^n$ induced by the injection $\bz^n\rightarrowtail \bz^d$ defined as 
$v_i\mapsto e_i, 1\leq i\leq n. $  This injection 
corresponds to an imbedding $\Gamma: G\rightarrowtail \wt{G}$.   
Identifying $\bz^d$ with $Hom(\bs^1,\wt{G})$, the splitting 
yields an identification  $\wt{G}\cong G\times H$,  
$\wt{g}=g.h$ where $g=\Gamma\circ\Lambda(\wt{g})\in G$ and $h=g^{-1}\wt{g}$.
The group $H\cong 
(\bs^1)^{d-n}$ is the subgroup 
of $\wt{G}$ with $\ker(\wt{\lambda})=Hom(\bs^1,H)\subset Hom(\bs^1,\wt{G})$. 
We let group $H$ act freely on the {\it right} of  
$\cz$ where $x.h=h^{-1}x\in \cz$ for $x\in \cz, ~h\in H\subset \wt{G}$. 
The quotient of $\cz$ by $H$ is the quasi-toric manifold $M$. 
(cf. \S4, \cite{dj}.)

Let $\chi:H\lr \bs^1$ be the restriction to $H$ of any character 
again denoted $\chi:\wt{G}\lr \bs^1$. 
One obtains a $G$-equivariant  
complex line bundle $L_\chi$ over $M$ 
with projection $\cz\times_H \bc_\chi\lr M$ where $\bc_\chi$ 
denotes the $1$-dimensional complex representation space 
corresponding to $\chi$.  Here the Borel construction $\cz\times_H\bc_\chi$
is obtained by the identification 
$$(xh,z)\sim(x,\chi(h)z), ~h\in H, ~x\in \cz, ~z\in \bc. \eqno(3.1)$$  

Equivalently $\cz\times_H\bc_\chi$ is the quotient 
of the diagonal action by $H$ on the left on 
$\cz\times \bc_\chi$.  The equivalence class of $(x,z)$ is denoted by 
$[x,z]$. 
The $G$-action on $L_\chi$ is given by 
$g.[x,z]:=[gx,\chi(g)z]$ for $x\in \cz,~z\in \bc_\chi$. 
     
When $\chi=\rho_j$ is the $j$-th projection $\wt{G}\lr \bs^1$,  
the corresponding $G$-line bundle on $M$  will be 
denoted $L_j$. Denote by $\pi_j:L_j\lr M$ the projection of the 
bundle $L_j$.

Henceforth we shall identify the character 
group of $\wt{G}$ with $Hom(\bz^d,\bz)$ etc. If $u\in Hom(\bz^d,\bz)$ 
vanishes on $\ker(\wt{\lambda})$, then the line bundle $L_u$ is  
isomorphic to the product bundle. However the $G$ action 
on it is given by the character $u|G$.  
 
Given $u\in Hom(\bz^n,\bz)=Hom(G,\bs^1)$, composing with surjection 
$\wt{G}\lr G$, we obtain a character of $\wt{G}$ which is trivial 
on $H$. As an element of $Hom(\bz^d,\bz)$, this is just the 
composition $\wt{u}=u\circ \wt{\lambda}$. 
Let  $e_1^*, 
\cdots, e_d^*$ be the dual of the standard basis for $\bz^d$. 
Note that the character $\wt{G}\lr \bs^1$ corresponding to $e_j^*$ 
is just the $j$-th projection $\rho_j$. Clearly,  
$$\wt{u}=u\circ\wt{ \lambda} 
=\sum_{1\leq j\leq d} u(\wt{\lambda}(e_j))e_j^*=
\sum_{1\leq j\leq d}u(v_j)e_j^*.$$
Hence we obtain the following isomorphism of $G$-bundles:
$$L_{\wt{u}}=\prod_{1\leq j\leq d}L_j^{u(v_j)}.\eqno(3.2)$$ 
Note that since $\wt{u}|H$ is trivial, $\cz\times_H\bc_{\wt{u}}
=M\times \bc$ and so $L_{\wt{u}}$ is isomorphic to the product bundle.    

Let $1\leq j\leq d.$ Choose an affine map  $h_j:\br^n\lr \br$
such that $h_j$ vanishes on $F_j$ and $h_j(p)>0$ for
$p\in P\setminus F_j$. Since  $\wt{G}_{F_j}:=\wt{G}_j$ acts freely 
on $\cz- p^{-1}(F_j)$, 
one has a well-defined trivialization 
$\sigma_j:\pi_j^{-1}(M-M_j)\lr(M-M_j)\times \bc_j$ given by 
$\sigma_j([x,z])=([x], \rho_j(g^{-1})z)$ where 
$x=[g,p]\in \cz, z\in \bc_j$. 

Using $\sigma_j$ and $h_j$ one obtains a well-defined section $s_j:M\lr L_j$ 
by setting $s_j([x])=[x,h_j(p)\rho_j(g)]$ where $x=[g,p]\in \cz$. 
Note that the section $s_j$ vanishes precisely on $M_j$. 
It is straightforward to verify that $s_j$ is $G$ equivariant. 

Now let $1\leq j_1,\cdots, j_k\leq d$ be such that $F_{j_1}\cap \cdots 
\cap F_{j_k}=\emptyset$. Thus $M_{j_1}\cap \cdots \cap M_{j_k}=\emptyset$. 
Consider the section $s:M\lr V$ defined as $s(m)=
(s_{j_1}(m),\cdots,s_{j_k}(m))$ where $V$ is (the total space of)  
the vector bundle $L_{j_1}\oplus \cdots \oplus L_{j_k}$.  
The section $s$ is nowhere vanishing: indeed $s(m)=0 \iff 
s_{j_r}(m)=0~\forall r \iff m\in M_{j_r}~\forall r$. Since 
$\cap_{1\leq r\leq k}M_{j_r}=\emptyset$, we see that $s$ is nowhere  
vanishing.  Since $V$ has geometric dimension at most $k-1$, 
applying the $\gamma^k$-operation to $[V]-k$ we obtain 
$\gamma^k([V]-k)=
\gamma^k(\oplus_{1\leq r\leq k}([L_{j_r}]-1))
=\prod_{1\leq r\leq k} ([L_{j_r}]-1)$. Therefore    
$$\prod_{1\leq r\leq k}(1-[L_{j_r}])=0\eqno(3.3)$$ 
whenever $\cap_{1\leq r\leq k}F_{j_r}=\emptyset$.

\noindent
\brem \label{normal}
(i) Let $\wt{L}_j$ denote the pull-back of $L_j$ by
the quotient map $\cz\lr M$. Since $H$ acts freely on
$\cz$, $\wt{L}_j$ is isomorphic to the product bundle.
This is the same as {\it dual} of the bundle $\wt{L}_j$ considered in
\S6.1 of \cite{dj}. A description of the stable tangent bundle of $M$ 
was obtained in Theorem 6.6 of \cite{dj}. It follows from their 
proof that the $L_j|M_j$ is isomorphic to 
the normal bundle to the imbedding $M_j\subset M$.  
Therefore we have $c_1(L_j)=e(L_j)=\pm[M_j]\in H^2(M;\bz)$  
where $e(L_j)$ denotes the Euler class of $L_j$ (see \cite{ms}.) 
The omni-orientation corresponding to $\lambda$ is so chosen 
as to have $c_1(L_j)=+[M_j]$.  

\noindent
(ii)  The complex projective $n$-space $\bp^n$ is a quasi-toric 
manifold over the standard $n$-simplex 
$\Delta^n=\{x=\sum_{1\leq i\leq n}x_ie_i\in \br^n\mid \sum_{1\leq i\leq n}
x_i\leq 1,~0\leq x_i\leq 1~~\forall i\geq 1\}$.  
The characteristic map $\lambda$ sends  
the facet  
$\Delta_i^n=\{x\in \Delta^n\mid x_i=0\} $ opposite the vertex 
$e_i, 1\leq i\leq n,$ to the 
standard basis element 
$v_i=e_i\in \bz^n$ for $i>0$ and sends the facet opposite the 
origin $\Delta_0^n=\{x\in \Delta^n|\sum_{1\leq j\leq n} x_j=1\}$ to the 
vector $v_0:=-(e_1+\cdots +e_n)\in \bz^n$. The line bundle 
$L_j$ is then verified to be isomorphic to the dual of 
tautological bundle over $\bp^n$ for $0\leq j\leq n$.  $M$ is canonically 
oriented as a complex manifold.  The omni-orientation corresponding 
to this choice of $\lambda$  
on $M$ determined by $\lambda$ is the orientation 
on $M_j\cong \bp^{n-1}$ determined by its  complex structure. 
(See also Example 5.19, \cite{bp}.) 
\erem 

\bpropo \label{kofm}
With notations as above, let $\cR=\cR(\bz;\lambda)$ where 
$r_i=1~\forall i\leq n$. One has 
a well-defined homomorphism $\psi:\cR\lr K(M)$
of rings which is in fact an isomorphism.
\epropo 
\noindent 
{\bf Proof:}  Relations (3.3) and (3.2) above clearly imply that 
$\psi$ is a well-defined algebra homomorphism. 

From Theorem 4.14 \cite{dj},  
the integral cohomology of $M$ is generated by degree $2$ elements.  
Indeed these can be taken to be dual cohomology classes 
$[M_j], 1\leq j\leq d$, where 
$M_j=\pi^{-1}(F_j)$, $F_j$ being the facets of $P$. 
As noted in Remark \ref{normal}, $c_1(L_j)=[M_j],~1\leq j\leq d.$    

Now Lemma 4.1, \cite{su}, implies that $K(M)$ is 
generated by the line elements $[L_j], 1\leq j\leq d$. 
This shows that $\psi$ is surjective. To show that it is 
injective, we observe that since $M$ is a CW complex 
with cells only in even dimensions, $K(X)$ is a free abelian 
group of rank $\chi(M)$ the Euler characteristic of $M$. 
But $\chi(M)=m,$ the number of vertices of $P$.  (In fact 
a $\bz$-basis for the integral cohomology of $M$ is 
the set of dual cohomology classes $[M_w], w\in P_0$.)  
Since by Lemma \ref{gen} the rank of $\cR$, as an 
abelian group, is {\it at most} $m$, it follows that 
$\psi$ is in fact an isomorphism of rings. \hfill $\Box$

We shall now prove the main theorem.

\noindent
{\bf Proof of Theorem \ref{main}:}  Note that the complex line
bundles $L_j$ are $G$-equivariant. Denote by $\cl_j$ the bundle
$E(L_j):=E\times_G L_j$ over $E(M)=E\times_G M$. Since the sections
$s_j$ are equivariant, so is the section $s=(s_{j_1}, \cdots, s_{j_k})$ of
$V=L_{j_1}\oplus\cdots\oplus L_{j_k}$. 
Hence we obtain a section $E(s):E(M)\lr E(V)$. If
$\bigcap_{1\leq r\leq k} F_{j_r}=\emptyset$, then $s$ and
hence $E(s)$ is nowhere vanishing. Again by applying the
$\gamma^k$-operation to $[E(V)]-k\in K(E(M))$, we conclude that 
$$\prod_{1\leq r\leq k}([\cl_{j_r}]-1)=0\eqno(3.4)$$
whenever $\bigcap_{1\leq r\leq k}F_{j_r}=\emptyset$. 

Since the isomorphism in equation (3.2) is
$G$-equivariant, one obtains an isomorphism
$\cl_u=\prod_{1\leq j\leq d}\cl_j^{u(v_j)}$.
Since $L_u$ is the product bundle $M\times \bc_u\lr M$,  
the bundle $\cl_u\cong p^*(\xi_1^{a_1}\cdots \xi_n^{a_n})$ 
where $a_i=u(v_i)$. (See \S1.) It follows that, in the 
$K(B)$-algebra $K(E(M))$ one has 
$$\prod_{1\leq j\leq d}[\cl_j]^{u(v_j)}=[\xi_1]^{a_1}\cdots [\xi_n]^{a_n}.
\eqno{(3.5)}$$ 

Setting $r_i=[\xi_i]$ for $1\leq i\leq n$, in view of (3.4) and (3.5) 
we see that there is a well-defined 
homomorphism of $K(B)$-algebras $\psi:\cR(K(B);\lambda) 
\lr K(E(M))$ which maps 
$x_j$ to $(1-[\cl_j])$ for $1\leq j\leq d$. 

It follows from Prop. \ref{kofm} that, the monomials in the $L_j$ generate 
$K(M)$.   Hence the fibre-inclusion $M\lr E(M)$ is totally non-cohomologous 
to zero in $K$-theory as the bundles $\cl_j$ restrict to $L_j$. 
As $B$ is compact Hausdorff and $K(M)$ is free abelian, we observe 
that the hypotheses of Theorem 1.3, Ch. IV, \cite{k}, are satisfied. 
It follows that $K(E(M))\cong K(B)\otimes K(M)$ as a $K(B)$-module. 
In particular, we conclude that $\psi$ is surjective. To see 
that $\psi$ is a monomorphism, note that $K(E(M))$ is a free 
module over $K(B)$ of rank $\chi(M)=m$, the number of vertices in $P$. 
Since, by Lemma \ref{gen} $\cR(K(B),\lambda)$ is generated by $m$ elements, 
it follows that $\psi$ is an isomorphism. \hfill $\Box$    

To conclude this section, we illustrate the above theorem 
in the case when $M$ is the $n$-dimensional complex projective 
space.  We remark that this special case follows immediately 
from Theorem 2.16, Ch. VI of \cite{k} as well.   

The complex projective $n$-space $\bp^n$ is a quasi-toric 
manifold over the standard $n$-simplex 
$\Delta^n=\{x=\sum_{1\leq i\leq n}x_ie_i\in \br^{n}\mid \sum_i x_i\leq 1,~
0\leq x_i\leq 1~~\forall i\geq 1\}$.  The characteristic map 
$\lambda$ sends  
the $i$-th facet --- the face opposite the vertex $e_i$ --- 
to the standard basis element  
$v_i:=e_i\in \bz^n$ for $i\geq 1$ and sends the $0$-th facet which is opposite 
the vertex $0$ to $v_0:=-(e_1+\cdots +e_n)\in \bz^n$.     
The space $E(\bp^n)$ is just the projective space bundle  
$\bp({\bf 1}\oplus \xi_1\oplus \cdots \oplus  
\xi_n)$ over $B$.  Here ${\bf 1}$ denotes the  
trivial complex line bundle $B\times \bc e_0\lr B$  
over $B$.  Indeed the map which sends 
$[e,x]=[(w_1,\cdots,w_n),[x_0,\cdots, x_n]]$ to the 
complex line spanned by the vector $x_0e_0+x_1w_1
+\cdots+x_nw_n$ in the fibre 
$\bp(\bc e_0+\bc w_1+\cdots +\bc w_n)$ over $\pi(e)\in B$, 
where $e=(w_1,\cdots, w_n)$, is a well defined bundle isomorphism.  
 
\beg \label{projbun} (Cf. Chapter IV, Theorem 2.16, \cite{k}.)  
$K(E(\bp^n))\cong K(B)[y]/\langle \prod_{0\leq i\leq n}(y-[\xi_i])
\rangle$ where $\xi_0:=1$ and $y$ is the class of the canonical 
bundle over $E(\bp^n)$ which restricts to the tautological 
bundle on each fibre of $E(\bp^n)\lr B$.  
\eeg 
\noindent 
{\bf Proof:}  
By choosing $u\in Hom(\bz^n,\bz)$ to be the 
dual basis element $e_i^*, i\geq 1,$ relation (1.3) gives 
$[\cl_i]=[\xi_i][\cl_0]$ in $\cR(K(B);\lambda)$ since $\xi_0=1$.    
It can be seen easily that  
other choices of $u$ in relation (1.3) do not lead 
to any new relation.  
Substituting for $[\cl_i]$ in 
relation (1.2),  we obtain 
$\prod_{0\leq i\leq n}(1-[\cl_0][\xi_i])=0$.  
Setting $y=[\cl_0]^{-1}$  we obtain 
$\prod([\xi_i]-y)=0$.  By theorem \ref{main} it 
follows that $K(E(\bp^n))\cong K(B)[y]/\langle \prod([\xi_i]-y)
\rangle.$  Since $y=[\cl_0^*]$, the proof is completed by observing that 
$\cl_0^*$ restricts to the {\it tautological} bundle on 
the fibres $\bp^n$. \hfill $\Box$

%%%%%%%%%%%%%%%%%%%%%%%%%%%%%%%%%%%%%%%%%%%%%%%%%%%%%%%%%%%%%%
%%%%%%%%%%%%%%%%%%%%%%%%%%%%%%%%%%%%%%%%%%%%%%%%%%%%%%%%%%%%%%

%%%%%%%%%%%%%%%%%%%%%%%%%%%%%%%%%%%%%%%%%%%%%%%%%%%%%%%%%%%%%%%%%%
%%%%%%%%%%%%%%%%%%%%%%%%%%%%%%%%%%%%%%%%%%%%%%%%%%%%%%%%%%%%%%%%%%
%%%%%%%%%%%%%%%%%%%%%%%%%%%%%%%%%%%%%%%%%%%%%%%%%%%%%%%%%%%%%%%%%%

\section{Bott-Samelson varieties}
In this section we illustrate our theorem in the case of
Bott-Samelson manifolds which were first 
constructed in \cite{bs} to study cohomology of 
generalized flag varieties. M.Demazure and D.Hansen used it to obtain 
desingularizations of Schubert varieties in 
generalized flag varieties. 
M.Grossberg and Y.Karshon \cite{gk} constructed Bott towers, which are 
iterated fibre bundles with fibre at each stage being $\bp_\bc^1$.  
They also showed that 
Bott-Samelson 
variety can be deformed to a toric variety.  The `special fibre,' of this 
deformation is a Bott tower. The underlying differentiable structure 
is preserved under the deformation. It follows that Bott-Samelson varieties 
have the structure of a   
quasi-toric manifold with quotient polytope being the $n$-dimensional 
cube $I^n$ where $n$ is the complex dimension of the Bott-Samelson variety.  
This quasi-toric structure has been  
explicitly worked out by Grossberg-Karshon \cite{gk}  and by 
M. Willems \cite{wil}. In this section we use  
Example \ref{projbun} to describe the $K$-ring of the Bott towers 
in terms of generators and relations.  Perhaps our theorem 
\ref{kbs} is well-known to experts but we could not find it 
explicitly stated in the literature. 

Let $C=(c_{i,j})$ denote an $n$-by-$n$ unipotent upper triangular matrix 
with integer entries.  The matrix $C$ determines 
a Bott tower $M(C)$ of (real) dimension 
$2n$.  Using the notation of 
\S3, it turns out that $\cz=(\bs^1)^{2n}\times I^{n}/\sim$ is 
the space $(\bs^3)^n\subset \bc^{2n}$. The quasi-toric manifold $M(C)$ 
is the quotient of $(\bs^3)^n$ by the 
action of $H=(\bs^1)^n$ on the right of $\cz$ where 
$$\begin{array}{lll}
(z_1, w_1,\cdots, z_n,w_n).t_i&=&(z_1,w_1,\cdots, z_it_i,w_it_i, 
\cdots z_j, w_jt_i^{c_{i,j}},\cdots, z_n,w_nt_i^{c_{i,n}})\cr
&=&(z_1,w_1,\cdots,z_it_i, w_it_i, \cdots, z_j,t_i^{-c_{i,j}}w_j,
\cdots, z_n,t^{-c_{i,n}}w_n)
\end{array} $$
for $t_i$ in the $i$-th factor of $H$.  The group $G=(\bs^1)^n$ acts 
on the left of $\cz=(\bs^3)^n$ by $t.(z_1,w_1,\cdots, z_n,w_n)
=(z_1,t_1^{-1}w_1, \cdots, z_n,t_n^{-1}w_n)$. This descends to 
an action on $M(C)$. 

For $1\leq i\leq n$ denote by $\cl_i$ the complex line 
bundle over $M(C)$ associated 
to the character $\rho_i:H\lr \bs^1$ defined as the  
projection to the $i$-th coordinate.  
Thus the total space of $\cl_i$ is obtained as the 
fibre product 
$\bs^3\times_H \bc$ where 
$((z_1,w_1,\cdots, z_n,w_n).t,\lambda)\sim((z_1,w_1,\cdots, z_n,w_n),t_i
\lambda)$  
for $t=(t_1,\cdots,t_n)\in H$ and $\lambda\in \bc$.  
We shall refer to $\cl_i$ as the $i$-th canonical bundle 
over $M(C)$. 
 
Let $C_i$ be the matrix obtained as the top $i$-by-$i$ diagonal 
block of $C,~ 1\leq i\leq n$.  Let $M(C_0)$ be the space consisting of 
a single point and let $M(C_1)=\bs^2=\bp^1_\bc$. 
We shall denote by $H_i=(\bs^1)^i$ the subgroup of $H$ where 
the last $n-i$ coordinates are the identity element.  
For each $i\geq 1$, one has the corresponding 
Bott tower $M(C_i).$  Let $L_j$ denote the $j$-th 
canonical bundle over $M(C_i)$, $1\leq j\leq i$. 
Consider the projection $\pi_i:M(C)\lr M(C_i)$ 
induced by the projection $(\bs^3)^n\lr (\bs^3)^i$ 
onto the first $i$-coordinates. Then $\pi_i^*(L_j)\cong \cl_j$ 
for $1\leq j\leq i$.  Indeed one has a commuting diagram   

$$\begin{array}{clc}
(\bs^3)^n\times_H\bc &\lr& (\bs^3)^i\times_{H_i}\bc\\
\downarrow&~&\downarrow\\
M(C)&\lr& M(C_i)
\end{array}$$  
where top horizontal map is the `bundle map' 
defined as $[(z_1,w_1,\cdots, z_n,w_n),\lambda]
\mapsto$ \linebreak $[(z_1,w_1,\cdots,z_i,w_i),\lambda]$, the vertical maps 
are projections of the bundle $\cl_j$ and $L_j$.   
It follows from 
Lemma 3.1, \cite{ms} that $\pi_i^*(L_j)\cong\cl_j.$ 
In view of this, by abuse of notation, we shall  
denote by the same symbol $\cl_j$ the 
$j$-th canonical bundle $L_j$ on $M(C_i).$   

Consider the map  
$\pi_{i,i+1}:M(C_{i+1})\lr M(C_{i})$ induced from the 
projection onto the first 
$i$-factors $(\bs^3)^{i+1}\lr (\bs^3)^i$. 
The map $\pi_{i,i+1}$ is the projection of the $\bs^2=\bp^1_\bc$-bundle 
associated with the complex vector bundle ${\bf 1}\oplus {\bf L}_i$ where 
${\bf L}_i$ is given by the character $(c_{1,i+1},\cdots, c_{i,i+1})
\in Hom(\bz^{i},\bz)=Hom(H^i,\bs^1)$ (see relation (3.1) in \S3).  
Thus 
$${\bf L}_i\cong \cl_1^{c_{1,i+1}}\cdots \cl_{i}^{c_{i,i+1}}.\eqno(4.1)$$
 
Denote by $\eta$ the complex line  bundle  
over $M(C_{i+1})$ which restricts to the dual of the tautological 
bundle on the fibres of $\pi_{i,i+1}:M(C_{i+1})\lr M(C_i)$.  
Observe that $\eta$ is just the bundle  
associated to the character $\rho_{i+1}:H^{i+1}\lr \bs^1$.  
Hence $\eta=\cl_{i+1}$ on $M(C_{i+1}).$

Let $i\geq 0$.  Note that one has $G$-equivariant sections 
$\sigma_i,\sigma'_i:M(C_i)\lr M(C_{i+1})$ of the bundle 
$M(C_{i+1})\lr M(C_i)$ 
defined as 
$\sigma_i([z_1,w_1,\cdots, z_i,w_i])=~[z_i,w_1,\cdots, z_i,w_i,1,0]$ and 
$\sigma'_i([z_1,w_1,\cdots, z_i,w_i])=[z_i,w_1,\cdots, z_i,w_i,0,1]$. 
The images of these sections are imbedded submanifolds of $M(C_{i+1})$ which 
correspond to the facets $I^{i}\times \{0\}$ and $I^i\times
\{1\}$ of $I^{i+1}$ respectively.  
We regard $M(C_i)$ as a submanifold of $M(C_{i+1})$ via $\sigma_i$.   
The normal bundle to the imbedding $M(C_i)\subset M(C_{i+1})$ is  
just the bundle $\cl_{i+1}|M(C_i)$.  (Cf. Remark 3.1.)
It follows that 
$$c_1(\cl_{i+1})=[M(C_i)]\in H^2(M(C_{i+1});\bz).\eqno(4.2)$$ 

For $1\leq i\leq n$ set $M_{i}:=\pi_i^{-1}(M(C_{i-1})).$  
Using the observations made above, the fact that $\pi_i$ is the projection 
of a fibre bundle, and equation  
(4.2), it follows that 
the line bundle $\cl_i|M_i$ is normal to the imbedding  
$M_{i}\subset M(C)$ and 
$$c_1(\cl_i)=[M_{i}]\in H^2(M(C);\bz). \eqno(4.3)$$

The structure of $M(C)$ as an 
iterated $2$-sphere bundle enables one to compute its $K$-ring 
using Example 3.3. 
   
\bth \label{kbt} We keep the above notations. 
Let $C$ be any $n\times n$ unipotent upper triangular matrix 
over $\bz$. The map $y_i\mapsto [\cl_i^*],~ 1\leq i\leq n$ defines 
an isomorphism of rings  from 
$K_n:=\bz[y_1^{\pm 1},\cdots, y_n^{\pm 1}]/J$ to $K(M(C))$  
where $J$ is the ideal generated by the elements: 
$(y_i-1)(y_i-y_0y_1^{-c_{1,i}}\cdots y_{i-1}^{-c_{i-1},i})$, $1\leq i\leq n,$
where $y_0:=1$.  One has $c_1(y_i)=-[M_{i-1}]\in H^2(M(C);\bz)$.
\eeth 
\noindent    
{\bf Proof:} This follows from Example \ref{projbun} by induction.  
Indeed, when $n=1$, $K(\bs^2)=K(\bp^1_\bc)\cong 
\bz[y_1]/\langle (y_1-1)^2\rangle \cong K_1$ since $y_1$ is invertible 
in $\bz[y_1]/\langle (y_1-1)^2\rangle$ where $y_1$ is the class of 
the dual of the tautological bundle over $\bp^1_\bc$. 
 
Let $i\geq 1$. By induction assume that  
$K_i\cong K(M(C_i))$ where $y_j\mapsto [\cl_j^*]$ for $1\leq j\leq i$ 
and that $y_j-1\in K_i$ is nilpotent for $j\leq i$. 
Since the $(i+1)$-st canonical bundle over $M(C_{i+1})$ is the 
bundle that restricts to the dual of the tautological bundle 
along the fibres $\bp_\bc^1$ of $\pi_{i,i+1}:M(C_{i+1})\lr M(C_i)$, 
from Example \ref{projbun} we obtain that $K(M(C_{i+1}))\cong K_{i}[y_{i+1}]/
\langle(y_{i+1}-1)(y_{i+1}-{\bf L}_i)\rangle$
under the $K_i$-algbera map that sends $y_{i+1}$ to $[\cl_{i+1}^*]$. 
Substituting for $[{\bf L}_i]$ from equation (4.1) 
%${\bf L}_i\cong\cl_1^{c_{1,i+1}}\cdots\cl_i^{c_{i,i+1}}$ by equation (4.1),  
we obtain the relation  
$(y_{i+1}-1)(y_{i+1}-y_1^{-c_{1,i+1}}\cdots y_i^{-c_{i,i+1}})$. 
Note that the (virtual) rank of $y_{i+1}-1$ is zero. It follows that 
$y_{i+1}-1$ is nilpotent.  One can 
also see this more directly by using induction and  
observing first that 
$y_1^{-c_{1,i+1}}\cdots y_{i}^{-c_{i,i+1}}-1$ can be expressed as a 
polynomial without constant term in $y_j-1$, $j\leq i$, 
and hence is nilpotent. Therefore  $y_{i+1}-1$ is also nilpotent.  
Hence $K_{i}[y_{i+1}]/\langle(y_{i+1}-1)
(y_{i+1}-y_1^{-c_{1,i+1}}\cdots y_{i}^{-c_{i,i+1}})\rangle=K_{i+1}$. Since 
$y_{i+1}-1$ has been shown to be nilpotent, the induction step is 
complete.  
The assertion about the first Chern class of $y_{i+1}$ is immediate from 
equation (4.3).  
\hfill $\Box$ 

Let $\cg$ be a complex semisimple linear algebraic group, $B$ a Borel 
subgroup. Fix an {\it algberaic} maximal 
torus $T\cong (\bc^*)^l$, $l$ being the rank of $\cg$, contained in 
$B$ and let $\Phi^+,\Delta$ denote 
the corresponding system of positive roots 
and simple roots respectively.  Denote by $W$ the Weyl group of $\cg$ with 
respect to $T$ and $S\subset W$ the set of simple reflections $s_\alpha,
\alpha\in \Delta$. 
For $\gamma\in \Delta$, denote 
by $P_\gamma\supset B$ the minimal parabolic subgroup  
corresponding to $\gamma$ so that $P_\gamma/B\cong \bp^1_\bc$.  
Let $\alpha_1,\cdots,\alpha_n$ be any sequence of simple roots. 
Consider the Bott-Samelson variety 
$M=P_{\alpha_1}\times_B\cdots\times_B P_{\alpha_n}\times_B\{pt\}$.  
Explicitly $M$ is the quotient of 
$\cp:=P_{\alpha_1}\times \cdots \times 
P_{\alpha_n}$ by the action of $B^n$ given by 
$(p_1,\cdots, p_n).b=(p_1b_1,b_1^{-1}p_2b_2, \cdots, b_{n-1}p_nb_n)$ 
for $(p_1,\cdots, p_n)\in \cp, 
(b_1,\cdots, b_n)\in B^n$.     
When $w=s_{\alpha_1}\cdots s_{\alpha_n}\in W$ is 
a reduced expression for $w$, one has a surjective birational morphism  
$M\lr X(w)$  which maps $[p_1,\cdots, p_n]$ to the coset 
$p_1\cdots p_n. B$ in the Schubert variety $X(w)\subset \cg/B$.  
In this case, $M$ is the Bott-Samelson-Demazure-Hansen 
\cite{bs},\cite{dem},\cite{han} desingularization of 
the Schubert variety associated to the reduced expression $w=s_{\alpha_1}
\cdots s_{\alpha_n}$. The 
Bott tower, which arises as the special fibre of a certain deformation of 
the complex structure of $M$,  
is associated to the unipotent upper triangular matrix $(c_{ij})$ where  
$c_{ij}=\langle \alpha_i,\alpha_j^\vee\rangle$, $i<j$.  The 
polytope which arises as the quotient of the Bott tower 
by the $n$-dimensional ({\it compact}) torus action is the $n$-cube $I^n$.  
See \cite{gk} or \cite{wil} for details.  
Feeding this data into Theorem \ref{kbt}, we 
obtain explicit description of the $K$-ring of a Bott-Samelson 
variety which is diffeomorphic to the Bott tower.  Alternatively one 
could use Example 3.3 and induction to obtain the same result. 

The Bott-Samelson variety  $M$ has an algebraic cell decomposition, i.e.,
a CW structure where the open cells are affine spaces contained in
$M$.  The closed cells of real codimension $2$ are the divisors $M_j$,
defined as the image of $\{(p_1,\cdots, p_n)\in \cp\mid p_j=1\}$ under
the canonical map $\cp\lr M.$  Note that the integral cohomology ring
of $M$ is generated by the dual cohomology classes $[M_j]\in
H^2(M;\bz)$, $1\leq j\leq n.$  

Lemma 4.2 of \cite{su} implies that the
forgetful map $\ck(M)\lr K(M)$ is an isomorphism of rings where
$\ck(M)$ is the Grothendieck $\ck$-ring of $M$.  
The following theorem is established using Theorem 4.1. 

\bth \label{kbs} 
Let $M$ be the (generalized) Bott-Samelson variety $P_{\alpha_1}\times_B 
\cdots\times_BP_{\alpha_n}\times_B\{pt\}$. Let $c_{i,j}=\langle 
\alpha_i,\alpha_j^\vee\rangle, 1\leq i<j\leq n$.  
The Grothendieck ring $\ck(M)$ of algebraic vector bundles on  
$M$ is isomorphic to $\bz[y_1^{\pm 1},\cdots, y_n^{\pm 1}]/
\langle (y_i-1)(y_i-y_0y_1^{-c_{1,i}}\cdots y_{i-1}^{-c_{i-1,i}}); 
1\leq i\leq n\rangle$ where $y_0:=1$.  The class $y_i$ is represented by the 
algebraic line bundle $\co(-M_i)$ for $1\leq i\leq n$. The forgetful 
ring homomorphism $\ck(M)\lr K(M)$ is an isomorphism.
\hfill $\Box$   
\eeth 

\brem
One has a well-defined involution $y_i\mapsto y_i^{-1}=:w_i$ of the algebra 
$K(M(C))$. Indeed multiplying the two factors in generating relation 
 $(y_i-1)(y_i-y_1^{-c_{1,i}}\cdots 
y_{i-1}^{-c_{i-1,i}})=0$ by $y_i^{-1}$ and $y_i^{-1}y_1^{c_{1,i}}\cdots 
y_{i-1}^{c_{i,i-1}}=0$ 
we get the same relation with the $y_j$'s replaced by $y_j^{-1}=w_j$: that is, 
$(w_i-1)(w_i-w_1^{-c_{1,i}}\cdots 
w_{i-1}^{-c_{i-1,i}})=0.$  Consequently, one could let $y_i$ to be 
the class of $\co(M_i)$ in the above theorem.  
\erem
\noindent 
{\bf Note:} After this work was completed, we come across
the papers of Civan and Ray \cite{cr} and M.Willems \cite{w2}.  
Civan and Ray determine the ring
structures of the generalized cohomology theories arising from complex
oriented ring spectra for Bott towers.  They also determine the
$KO$-rings of Bott towers using entirely different methods.  
Willems \cite{w2}  has computed the equivariant K-ring of Bott towers 
and Bott-Samelson varieties. While we
consider only $K$-ring, our results apply to the more general class of
quasi-toric manifolds.

\noindent
{\bf Acknowledgements:} We thank Prof. V.Balaji for 
valuable discussions.  We thank Prof. M.Masuda for his valuable comments 
and for careful reading of an earlier version of this paper.  
  
%%%%%%%%%%%%%%%%%%%%%%%%%%%%%%%%%%%%%%%%%%%%%%%%%%%%%%%%%%%%%%%%%%%%%%%%%%
%%%%%%%%%%%%%%%%%%%% REFERENCES %%%%%%%%%%%%%%%%%%%%%%%%%%%
%%%%%%%%%%%%%%%%%%%%%%%%%%%%%%%%%%%%%%%%%%%%%%%%%%%%%%%%%%%%%%%%%%%%%%%%%%% 

\noindent
Institute of Mathematical Sciences\\
CIT Campus, Chennai 600 113, INDIA \\
E-mail: {\tt sankaran@imsc.res.in}\\
        {\tt uma@imsc.res.in}\\

\noindent
Current Address (V.U.):\\
{ Institut Fourier\\
100, rue des Maths, BP74\\
38402 St Martin d'Heres\\
Cedex, FRANCE\\
E-mail:{\tt uma@mozart.ujf-grenoble.fr}} 

\begin{thebibliography}{9}
\bibitem {bb} A.Bahri and M.Bendersky, The $KO$-theory of toric 
manifolds, Trans. Amer. Math. Soc. {\bf 352}, 1191-1202.    
\bibitem{b} A.Bronsted, An introduction to convex polytope, (1983), 
Springer-Verlag, NY.
\bibitem{bp} V.M.Buchstaber and T.E.Panov, {\it Torus actions and 
their applications in topology and combinatorics}, 
Univ. Lect. Series-{\bf 24},(2002), AMS, Providence, RI.
\bibitem{bs} R.Bott and H.Samelson, Applications of the theory of Morse 
to symmetric spaces.  Amer. J. Math. {\bf 80},  1958, 964--1029. 
\bibitem{c} Y.Civan, Some examples in toric geometry, 
arXiv:math.AT/0306029 v2. 
\bibitem{cr} Y.Civan and N.Ray, Homotopy decompositions and $K$-theory 
of Bott towers, arXiv:math.AT/0408261 v1.
\bibitem{dj} M. W. Davis and T. Januszkiewicz, Convex polytopes, Coxeter 
orbifolds and torus actions, Duke Math. Jour. {\bf 62},(1991), 417-451.
\bibitem{dem} M.Demazure, 
D\'esingularisation des vari\'et\'es de Schubert g\'en\'eralis\'es,
Ann. Sci. \'Ecole Norm. Sup. {\bf 7} (1974), 53--88.

\bibitem{f} W.Fulton, Introduction to toric varieties, Ann Math Studies {\bf 
131},(1993), Princeton Univ. Press, Princeton, NJ. 
\bibitem{gk} M. Grossberg and Y. Karshon, Bott towers, 
complete integrability, and the extended character of representations,  
Duke Math. Jour.  {\bf 76},  (1994), 23--58. 
\bibitem{han} H.C.Hansen, On cycles in flag manifolds,  Math. Scand. {\bf 33}  
(1973), 269--274 (1974). 
\bibitem{k} M.Karoubi, {\it K-Theory}, 
Grundlehren der Mathematischen Wissenschaften {\bf 226}, 
Springer-Verlag, Berlin, 1978. 
%\bibitem{kk} B.Kostant and S.Kumar, $T$-Equivariant $K$-theory 
%of generalized flag varieties, J.Diff.Geom., {\bf 32}, (1990), 549-603.  
\bibitem{kh} A.G. Khovanskii, Hyperplane sections of polyhedra, 
Funct. Analys. Appl. {\bf 20}, (1986), 41-50. 
\bibitem{ms} J.W.Milnor, J.D.Shasheff, Characteristic classes, 
Ann. Math. Studies {\bf 76},(1974), Princeton 
Univ. Press, Princeton, NJ. 
\bibitem{su} P.Sankaran and V.Uma, Cohomology of toric bundles, 
Comment. Math. Helv., {\bf 78},(2003), 540-554. Errata, {\bf 79},(2004), 
840-841.
%\bibitem{seg} G.Segal, Equivariant $K$-theory, Publ. Math.,IHES, 
%{\bf 34}, (1968), 129-151.
\bibitem{wil} M.Willems, Cohomologie et K-th\'eorie \'equivariantes 
des tours de Bott et des vari\'et\'es de drapeaux. Application au 
calcul de Schubert, arXiv:math.AG/0311079 v1.
\bibitem{w2} M.Willems, K-theorie equivariante des varietes de Bott-Samelson.
Application a la structure multiplicative de la K-theorie equivariante 
des varietes de drapeaux, arXiv:math.AG/0412152.  
\bibitem{z} G.Ziegler, {\it Lectures on Polytopes,}  
Graduate Texts in Mathematics, 152. Springer-Verlag, New York, 1995. 
\end{thebibliography}
\end{document}